\documentclass[3p,times]{elsarticle_preprint}

\usepackage{amsfonts,amsmath,amssymb,amsthm}
\usepackage{graphicx}
\usepackage{verbdef}
\usepackage{algorithm}
\usepackage{algorithmic}
\usepackage{hyperref}

\newcommand{\RR}{\mathbb{R}}

\usepackage{color}

\newcommand{\bs}[1]{\boldsymbol{#1}}
\newcommand{\ns}{\mathcal N_{\Phi}(\Omega)}

\newtheorem{definition}{Definition}[section]
\newtheorem{property}{Property}[section]
\newtheorem{theorem}{Theorem}[section]

\newtheorem{remark}{Remark}[section]

\journal{R. Cavoretto, S. De Marchi, A. De Rossi, E. Perracchione, G. Santin}

\usepackage{amssymb}

\usepackage[figuresright]{rotating}

\begin{document}
\begin{frontmatter}

\title{Partition of unity interpolation using stable kernel-based techniques}


\author[label1]{R. Cavoretto\corref{cor1}}
\ead{roberto.cavoretto@unito.it}
\cortext[cor1]{Corresponding author.}

\author[label2]{S. De Marchi}
\ead{demarchi@math.unipd.it}

\author[label1]{A. De Rossi}
\ead{alessandra.derossi@unito.it}

\author[label1]{E. Perracchione}
\ead{emma.perracchione@unito.it}

\author[label3]{G. Santin}
\ead{santinge@mathematik.uni-stuttgart.de}

\address[label1]{Department of Mathematics \lq\lq G. Peano\rq\rq, University of Torino, via Carlo Alberto 10, I--10123 Torino, Italy}

\address[label2]{Department of Mathematics, University of Padova, via Trieste 63, I--35121 Padova, Italy}

\address[label3]{Institute of Applied Analysis and Numerical Simulation, University of Stuttgart, Pfaffenwaldring 57, D--70569 Stuttgart, Germany
	\vskip 0.3cm	
	\large{This paper is dedicated to Prof. Francesco A. Costabile on the occasion of his 70th birthday}}

\begin{abstract}
In this paper we propose a new stable and accurate approximation technique which is extremely effective for interpolating large scattered data sets. The Partition of Unity (PU) method is performed considering Radial Basis Functions (RBFs) as local approximants and using locally supported weights. In particular, the approach consists in computing, for each PU subdomain, a stable basis. Such technique, taking advantage of the local scheme, leads to a significant benefit in terms of stability, especially for flat kernels. Furthermore, an optimized searching procedure is applied to build the local stable bases, thus rendering the method more efficient. 
\end{abstract}


\end{frontmatter}




\section{Introduction} \label{intro}
Considering the state of the art \cite{Cavoretto14b,Cavoretto15b,Demarchi13,Demarchi15,Pazouki11}, we propose a new method for multivariate approximation which allows to interpolate large scattered data sets stably, accurately and with a relatively low computational cost.

The interpolant we consider is expressed as a linear combination of some basis or kernel functions. Focusing on Radial Basis Functions (RBFs), the Partition of Unity (PU) is performed by blending RBFs as local approximants and using locally supported weight functions. With this approach a large problem is decomposed into many small problems, \cite{Babuska97,Fasshauer07,Melenk96,Wendland02a}, and therefore in the approximation process we could work with a large number of nodes.
	
However, in some cases, local approximants and consequently also the global one suffer from instability due to ill-conditioning of the interpolation matrices. This is directly connected to the order of smoothness of the basis function and to the node distribution. It is well-known that the stability depends on the flatness of the RBF.  More specifically, if one keeps the number of nodes fixed and  considers smooth basis functions, then the problem of instability  becomes evident for small values of the shape parameter. Of course, a basis function with a finite order of smoothness can be used to improve the conditioning but the accuracy of the fit gets worse. For this reason, the recent research is moved to the study of more stable bases. 

For particular RBFs, techniques allowing to stably and accurately compute the interpolant, also in the \textit{flat limit} $\varepsilon\to 0$, have been designed in the recent years. These algorithms, named RBF-QR methods, are all rooted in a particular decomposition of the kernel, and they have been developed so far to treat the Gaussian and the {Inverse MultiQuadric kernel}. We refer to \cite{Fasshauer12,Fornberg04,Fornberg07,Fornberg11} for further details on these methods.

A different and more general approach, consisting in computing, via a truncated Singular Value Decomposition (SVD) stable bases, namely  Weighted SVD (WSVD) bases, has been presented in \cite{Demarchi13}. 
We remark that in the cases where the RBF-QR algorithms can be applied, they produce a far more stable solution of the interpolation problem. Nevertheless, the present technique applies to \textit{any} RBF kernel, and to any domain.


In this paper, a stable approach via the PU method, named WSVD-PU, which makes use of local WSVD bases and uses compactly supported weight functions, is presented. Thus, following \cite{Demarchi15}, for each PU subdomain a stable RBF basis is computed in order to solve the local interpolation problem. Consequently, since the local approximation order is preserved for the global fit, the interpolant results more stable and accurate. Concerning the stability, we  surely expect a more significant improvement in the stabilization process with infinitely smooth functions  than with functions characterized by a finite order of regularity. Moreover, in terms of accuracy, the benefits coming from the use of such stable bases are more significant in a local approach than in a global one. In fact, generally, while in the global case a large number of truncated terms of the SVD must be dropped to preserve stability, a local technique requires only few terms are eliminated, thus enabling the method to be much more accurate.

Concerning the computational complexity of the algorithm, the use of the so-called block-based space partitioning data structure enables us to efficiently organize points among the different subdomains, \cite{Cavoretto15c}. Then, for each subdomain a local RBF problem is solved with the use of a stable basis. The main and truly high cost, involved in this step, is the computation of the SVD. To avoid this drawback, techniques based on Krylov space methods are employed, since they turn out to be really effective, \cite{Beatson99,Demarchi15}. A complexity analysis supports our findings.

The guidelines of the paper are as follows. In Section \ref{WSVD}, we present the WSVD bases, computed by means of the  Lanczos algorithm, in the general context of global approximation. Such method is used coupled with the PU approach which makes use of an optimized searching procedure, as shown in Section \ref{PUM_WSVD}. The proposed approach turns out to be stable and efficient, as stressed in Section \ref{compl_anal_WSVD}. In Sections \ref{NE_LAN} and \ref{applicazione} extensive numerical experiments and applications, carried out with both globally and compactly supported RBFs of different orders of smoothness, support our results. 
Moreover, all the \textsc{Matlab} codes are made available to the scientific community in a downloadable free software package: 
\begin{center}
http://hdl.handle.net/2318/1527447.
\end{center}

\section{RBF interpolation and WSVD basis} \label{WSVD}

In Subsection \ref{rbf_prelim} we briefly review the main theoretical aspects concerning RBF interpolation, \cite{Buhmann03}, while the remaining subsections are devoted to the efficient computation of the WSVD basis via Krylov space methods.

\subsection{RBF interpolation}
\label{rbf_prelim}

Our goal is to recover a function $f: \Omega\to \RR$, $\Omega$ being a bounded set in $\RR^M$, using a set of samples of $f$ on $N$ pairwise distinct points $X_N\subset\Omega$, namely $\boldsymbol{f} =[f_1, \ldots , f_N]^T$, $f_i = f(\boldsymbol{x}_i)$, $\boldsymbol{x}_i\in X_N$. To this end, one considers a positive definite and symmetric kernel $\Phi : \Omega \times \Omega \longrightarrow \RR$ 
to construct an interpolant in the form 
\begin{equation}\label{ansatz}
R^N(\bs{x}) = \sum_{j = 1}^N c_j \Phi(\bs{x}, \bs{x}_j), \;\; \bs{x}\in\Omega.
\end{equation}
The kernels we consider are always radial, meaning that there exist a positive \textit{shape parameter} $\varepsilon$ and a function $ \phi: \RR_{\geq 0}\to \RR$ such that for all $\boldsymbol{x},\boldsymbol{y} \in \Omega$, $\Phi(\boldsymbol{x},\boldsymbol{y})=\phi_{\varepsilon}(||\boldsymbol{x}-\boldsymbol{y}||_2)=\phi(\varepsilon ||\boldsymbol{x}-\boldsymbol{y}||_2)$. In Table \ref{tab_rbf} we report a list of some strictly positive definite radial kernels with their smoothness degrees. We remark that Gaussian, Inverse MultiQuadric and Mat$\acute{\text{e}}$rn functions are globally supported and strictly positive definite in $\RR^M$ for any $M$, whereas Wendland ones are compactly supported (whose support is $\left[0,1/\varepsilon\right]$) and strictly positive definite in $\RR^M$ for $M\leq 3$ (see \cite{Wendland05}).

\begin{table}[ht!]
\begin{center}
\begin{tabular}{ll}
\hline
\rule[0mm]{0mm}{3ex}
RBF  & $\phi_{\varepsilon}(r)$ \\
\hline
\rule[0mm]{0mm}{3ex}
{\rm Gaussian $C^{\infty}$} (GA) & ${\rm e}^{-\varepsilon^2 r^2}$   \\
\rule[0mm]{0mm}{3ex}
{\rm Inverse MultiQuadric $C^{\infty}$} (IMQ) & $(1+\varepsilon^2r^2)^{-1/2}$   \\
\rule[0mm]{0mm}{3ex}
{\rm Mat$\acute{\text{e}}$rn $C^6$} (M6)  & ${\rm e}^{-\varepsilon r} (\varepsilon^3 r^3 + 6\varepsilon^2r^2+15\varepsilon r+15)$   \\
\rule[0mm]{0mm}{3ex}
{\rm Mat$\acute{\text{e}}$rn $C^4$} (M4)  & ${\rm e}^{-\varepsilon r} (\varepsilon^2r^2+3\varepsilon r+3)$   \\
\rule[0mm]{0mm}{3ex}
{\rm Wendland $C^6$} (W6)  & $\left(1-\varepsilon r\right)_+^8(32\varepsilon^3r^3+25\varepsilon^2r^2+8\varepsilon r+1)$  \\
\rule[0mm]{0mm}{3ex}
{\rm Wendland $C^4$} (W4)  & $\left(1-\varepsilon r\right)_+^6(35\varepsilon^2 r^2+18\varepsilon r+3)$  \\
\hline
\end{tabular}
\end{center}
\caption{Examples of strictly positive definite radial kernels with their orders of smoothness and shape parameter $\varepsilon > 0$; $r=||\cdot||_2$ is the Euclidean distance, while $(\cdot)_+$ denotes the truncated power function.}
\label{tab_rbf}
\end{table}

The real coefficients $\boldsymbol{c}= [c_1, \ldots, c_N]^T$ in \eqref{ansatz} are determined by solving the linear system $A \boldsymbol{c}= \boldsymbol{f}$, where the interpolation (or kernel) matrix $A \in \RR^{N \times N}$ is given by
\begin{equation}\label{linear system}
A_{ij}= \Phi (\boldsymbol{x}_i , \boldsymbol{x}_j), \quad i,j=1, \ldots, N.
\end{equation}
{The so constructed solution $R^N$ is a function of the \textit{native Hilbert space} ${\cal N}_{\Phi}(\Omega)$ uniquely associated with the kernel, and, if $f\in\ns$, it is in particular the $\ns$-projection of $f$ into the subspace 
	$$
	{\cal N}_{\Phi}(X_N)= \textrm{span} \{\Phi(\boldsymbol{x},\boldsymbol{x}_j), \boldsymbol{x}_j \in X_N\}.
	$$
}

Although this interpolation method is known to be highly unstable in most cases being the matrix $A$ severely ill conditioned, it has been proven (see \cite{DeMSch2010}) that the interpolation operator $f\mapsto R^N$ is stable as an operator in the function space $\ns$. This gap has been widely recognized to be caused by the use of the standard basis, and a lot of efforts have been made in recent years to introduce better or perfectly conditioned basis (see \cite{Pazouki11} for a general theoretical treatment of this topic, and \cite{Cavoretto15a, Fasshauer12, Fornberg04, Fornberg07} for particular instances of stable basis; for an overview see the book \cite{Fasshauer15}). 

\subsection{WSVD basis}

We are interested here in the use of the \textit{WSVD basis} introduced in \cite{Demarchi13}, thanks to its flexibility with respect to the choice of the kernel $\Phi$. We recall in the following some relevant properties of this basis, while we refer to the paper \cite{Demarchi13} for further details.

To construct a basis $\mathcal U = \{u_j\}_{j=1}^{N}$ of ${\cal N}_{\Phi}(X_N)$ it is enough to assign an invertible coefficient matrix $D_{\mathcal U} = [d_{ij}]_{i,j=1}^N$ such that 
$$
u_j(\boldsymbol{x}) = \sum_{i=1}^N d_{ij} \Phi(\boldsymbol{x}, \boldsymbol{x}_i), 
$$
or, equivalently, an invertible value matrix $V_{\mathcal U} = [u_{j}(\boldsymbol{x}_i)]_{i,j=1}^N$. The two matrices are related as $A = V_{\mathcal U} \cdot D_{\mathcal U}^{-1}$ (see \cite{Pazouki11}), and in our situation they are defined as follows (see \cite{Demarchi13}). 
\begin{definition}
	A WSVD basis ${\cal U}$ is a basis for ${\cal N}_\Phi( X_N)$ characterized by the matrices
	$$
	D_{\cal U} = \sqrt{{\tilde W}} \cdot Q \cdot\Sigma^{-1/2}\; \mbox{ and }\; V_{\cal U} = \sqrt{{\tilde W}^{-1}} \cdot Q \cdot \Sigma^{1/2}, \quad \textrm{where} \quad 	\sqrt{{\tilde W}} \cdot A \cdot \sqrt{{\tilde W}} =  Q \cdot \Sigma \cdot Q^{T},
	$$	
	is a singular value decomposition of the scaled kernel matrix $A_{{\tilde W}}=\sqrt{{\tilde W}} \cdot A \cdot \sqrt{{\tilde W}}$, and ${\tilde W}_{ij} = \delta_{ij} {\tilde w}_i$ is a diagonal matrix of positive weights. 
\end{definition}
We observe that the definition uses a set of positive weights that was employed in the original formulation as cubature weights to construct the basis. Nevertheless, these weights do not change the numerical behavior of the basis, hence we assume from now on ${\tilde w}_i = 1/N$, $1\leq i\leq N$. Moreover, for notational convenience, the diagonal elements of $\Sigma$ are denoted as $\sigma_1\geq\dots\geq\sigma_N$.

This basis has been introduced to mimic in a discrete sense the properties of the \textit{eigenbasis}. Such basis is constructed starting from the operator $T: L_2(\Omega) \to L_2(\Omega)$,
\begin{equation}\label{cont_vers}
	T[f](\boldsymbol{x}) = \int_{\Omega} \Phi(\boldsymbol{x},\boldsymbol{y})  f(\boldsymbol{y})  d\boldsymbol{y}, 
\end{equation}
through the following Theorem (see e.g. \cite{Pogo66}). 

\begin{theorem}[Mercer's Theorem]
If the kernel $\Phi$ is continuous and positive definite on a bounded set $\Omega \subseteq \mathbb{R}^{M}$, the operator $T$ 
has a countable set of eigenfunctions $\{ \varphi_k \}_{k}$ and eigenvalues	$\{\lambda_k \}_{k}$. The eigenfunctions are orthonormal in $L_2(\Omega)$ and orthogonal in $\ns$ with 
{$\|\varphi_k\|^2_{\ns} = \lambda_k^{-1}$.}
Moreover, the kernel can be expressed in terms of the eigencouples as
$$
\Phi(\boldsymbol{x}, \boldsymbol{y}) = \sum_{k}\lambda_k \varphi_k(\boldsymbol{x})\varphi_k(\boldsymbol{y}), \;\boldsymbol{x}, \boldsymbol{y}\in\Omega,
$$
where the series converges uniformly and absolutely. 
\end{theorem}
For its use in interpolation in $\ns$, it is convenient to use the basis $\{\sqrt{\lambda_k}\varphi_k\}_k$, that is normalized in $\ns$. With this normalization, the basis has the following properties.
\begin{property}
The eigenbasis $\{\sqrt{\lambda_k}\varphi_k\}_k$ has the following properties:
\begin{enumerate}[i.]
\item  it is $\ns$-orthonormal,
\item  it is $L_2(\Omega)$-orthogonal with {squared} norm $\lambda_k$,
\item $(\sqrt{\lambda_k}\varphi_k, f)_{L_2(\Omega)} = \lambda_k(\sqrt{\lambda_k}\varphi_k, f)_{\ns}$, $\forall f \in\ns$,
\item $\lambda_k\geq\lambda_{k+1}$ and $\lambda_k\to 0$ as $k\to\infty$,
\item $\sum_k\lambda_k = \phi(0)\ \mathrm{meas}(\Omega)$.
\end{enumerate}
\end{property}
As proven in \cite{Demarchi13}, the WSVD basis enjoys the same properties when the inner product of $L_2(\Omega)$ is replaced with its discrete version $\ell_2(X_N)$, as summarized in the following statement.
\begin{property}\label{discrete properties}
The WSVD basis $\{u_k\}_{k=1}^N$ has the following properties:
\begin{enumerate}[i.]
\item  it is $\ns$-orthonormal,
\item  it is $\ell_2(X_N)$-orthogonal with norm $\sigma_k$,
\item\label{double} $(u_k, f)_{\ell_2(X_N)} = \sigma_k(u_k, f)_{\ns}$, $\forall f \in\ns$,
\item $\sigma_1\geq\dots\geq \sigma_N>0$,
\item\label{discrete sum condition} $\sum_{k=1}^{N}\sigma_k = \phi(0)\ \mathrm{meas}(\Omega)$.
\end{enumerate}
\end{property}
Since the interpolation is a $\ns$-projection, we can rewrite the interpolant $R^N$ in terms of the $\ns$-orthonormal WSVD basis as
\begin{equation}\label{interp}
R^N(\boldsymbol{x})= \sum_{k=1}^{N} (f, u_k)_{\ns} u_k(\boldsymbol{x}),
\end{equation}
and, thanks to point \eqref{double} of Property \ref{discrete properties}, this can be further rewritten as
\begin{equation}
R^N(\boldsymbol{x})=\sum_{k=1}^{N} \sigma^{-1}_k (f, u_k)_{\ell_2(X_N)} u_k(\boldsymbol{x}).
\label{st_app}
\end{equation}
The latter form of the interpolant shows that $R^N$ is also the solution of the discrete least-squares approximation problem 
$\min|| f-g||_{\ell_2(X_N)}$
among all functions in $g\in\mathcal N_{\Phi}(X_N)$. If we instead solve the minimization problem over the subspace $\textrm{span} \{u_1, \ldots, u_m \}$, $m\leq N$, we find a solution $R^m$ given by the truncation of the interpolant, i.e., 
\begin{equation}\label{interp1}
R^m(\boldsymbol{x})=\sum_{k=1}^{m} \sigma^{-1}_k (f, u_k)_{\ell_2(X_N)} u_k(\boldsymbol{x}).
\end{equation}
Observe that it makes sense to consider the last minimization problem and its solution $R^m$ for some $m\leq N$ instead of the original interpolation problem, since in this way we leave out the portion of the subspace $\mathcal N_{\Phi}(X_N)$ corresponding to small singular values, and this corresponds to solve the linear system \eqref{linear system} by means of
a low-rank approximation of the matrix $A$. A detailed discussion of this approach can be found in \cite{Demarchi13}, but we remark here that this method strictly depends on the behavior of the singular values of $A$. Namely, if we consider smoother RBFs, then by using \eqref{interp1} instead of the standard approach a better stabilization of the interpolation process is expected. 

On the other hand, this method has some disadvantages. First, it is required to compute a singular value decomposition of the (possibly large) kernel matrix, and in the end only a few elements of the decomposition are used. This is computationally expensive, but in the next section we explain how to overcome this problem. Second, this method requires to neglect part of the information to reduce instability, and, in some cases, this removal is too big to obtain a meaningful approximant. A solution to this problem is provided by the coupling with a localization method, and it is the main topic of this paper.

\subsection{Fast computation through Krylov space methods}\label{WSVD_LAN}

We present here a way to compute an approximation of the WSVD basis that makes use of the Lanczos algorithm. The method is discussed in \cite{Demarchi15}, and it aims at reducing the computational cost of the procedure by approximating the truncated SVD of $A$.

We start by a general description of the Lanczos algorithm. Further details can be found in \cite{Chen09,Simon00,Wei10}. Let ${\cal K}_m(A,\boldsymbol{f}) = \textrm{span} \{\boldsymbol{f}, A\boldsymbol{f},\ldots, A^{m-1} \boldsymbol{f}\}$ be the Krylov subspace of order $m$ generated by the matrix $A$ and the vector $\boldsymbol{f}$. The Lanczos method computes an orthonormal basis $ \{ \boldsymbol{p}_i \}^m_{i=1}$ of  ${\cal K}_m(A,\boldsymbol{f})$ through a Gram-Schmidt orthonormalization, i.e., the Lanczos basis $\{ \boldsymbol{p}_i \}^m_{i=1}$, is computed by the following recurrence formula:
\begin{equation}
\beta_{i+1}\boldsymbol{p}_{i+1}=A\boldsymbol{p}_{i}-\alpha_{i}\boldsymbol{p}_i-\beta_{i}\boldsymbol{p}_{i-1}, \quad \textrm{with} \quad \beta_{1}\boldsymbol{p}_{0}=0. 
\end{equation}

Letting $P_m \in \mathbb{R}^{N \times m}$ the matrix having the vectors $\boldsymbol{p}_i$ as columns, and letting $H_m$ be the $(m+1)\times m$ tridiagonal matrix defined as
\begin{equation}
H_{m} =
\begin{pmatrix}
\alpha_1 & \beta_2 & \cdots & 0 \\
\beta_2 & \alpha_2 & \cdots & 0 \\
\vdots  & \vdots  & \ddots & \beta_{m-1}\\
0 &  \cdots &  \beta_{m-1} & \alpha_m\\
0 & \cdots & 0 & \beta_m
\end{pmatrix}
\end{equation}
the algorithm can be formulated in matrix form as 
\begin{equation}
AP_m = P_{m+1} \bar{H}_m, \quad \bar{H}_m = \begin{pmatrix}
H_m\\
\bar{h}\boldsymbol{e}^{T}_m
\end{pmatrix},
\end{equation}
where $\boldsymbol{e}_m \in \mathbb{R}^m$ is the unit vector and $\bar{h}$  is a scalar value.

Once we compute the matrices, the solution of the initial system can be approximated as $\boldsymbol{x} = P_m \boldsymbol{y}$, where $\boldsymbol{y} \in \mathbb{R}^m$ is such that $\bar{H}_m \boldsymbol{y}= ||\bs{f}||_2 \boldsymbol{e}_1$.

The idea is to use the matrix $\bar{H}_m$ to approximate the SVD of $A$, and since $A$ has usually a good low-rank approximation, we expect to do so with $m \ll N$. Specifically, let $\bar{H}_m=U_m \bar{\Sigma}_m V^T_m$ be a singular value decomposition of $\bar{H}_m$, where $U_m \in \mathbb{R}^{(m+1) \times (m+1)}$, $V_m \in \mathbb{R}^{m \times m}$ are unitary matrices and 
\begin{displaymath}
\bar{\Sigma}_m= \begin{pmatrix}
\Sigma_m\\
0\\
\end{pmatrix},
\end{displaymath}
with $\Sigma_m$ the diagonal matrix with singular values $\sigma_i$, $1\leq i\leq m$ on the diagonal.

Since the last row of $\bar{\Sigma}_m$ is the zero vector, the decomposition does not change if we remove this row and the last column of $U_m$.  Thus, to simplify the notation we denote by $U_m$ the matrix without the last column, so that the decomposition becomes $\bar{H}_m=U_m \Sigma_m V^T_m$.

Now we want to define an approximation of the WSVD basis using the approximate SVD of $A$. We define a set of functions $\{{\bar{u}}_k\}^m_{k=1} \in {\cal N}_{\Phi}(X_N)$ which shows similarities with the WSVD basis, even if  it does not form a basis, since they do not span $ {\cal N}_{\Phi}(X_{N})$. {Anyway, with abuse of notation, we call again this set of functions a \emph{basis}.}

\begin{definition}
	The approximate WSVD basis is characterized by the matrices
	$$
	D_{\bar{\cal U}_m}=P_m \cdot V_m \cdot \Sigma^{-1/2}_m\; \mbox{ and }\; V_{\bar{\cal U}_m}=P_{m+1}  \cdot U_m \cdot \Sigma_m^{1/2}.
	$$
	where $A P_m = P_{m+1} \bar{H}_m$ is the Lanczos decomposition of $A$ of order $m$ and $\bar{H}_m = U_m \Sigma_m V_m^T$ is a singular value decomposition of $\bar{H}_m$.	
\end{definition}
The next statement clarifies the connection between the two basis. As it is evident by construction, the basis strongly depends on the particular function $f \in {\cal N}_{\Phi}(\Omega)$. 

\begin{property}\label{approx discrete properties}
The approximate WSVD basis $\{\bar u_k\}_{k=1}^m$ has the following properties:
\begin{enumerate}[i.]
\item  it is near $\ns$-orthonormal, meaning that its $\ns$-Gramian is the identity matrix plus a rank one matrix,
\item  it is $\ell_2(X_N)$-orthogonal with norm $\sigma_k$,
\item\label{approx double} $({\bar u}_k, f)_{\ell_2(X_N)} = \sigma_k({\bar u}_k, f)_{\ns}$ if $f$ is the function used to construct the basis,
\item $\sigma_1\geq\dots\geq \sigma_m>0$,
\item\label{m=N} it coincides with the WSVD basis if $m = N$.
\end{enumerate}
\end{property}
This basis allows to solve again the least square approximation problem. Namely, if $f\in\ns$ is the function used for the Lanczos algorithm, the approximant $\bar{R}^m$ defined as
\begin{equation*}
\bar{R}^{m}(\boldsymbol{x})=\sum_{k=1}^{m} \sigma^{-1}_k (f, \bar{u}_k)_{\ell_2(X_N)} \bar{u}_k(\boldsymbol{x}), 
\end{equation*}
minimizes the distance $\|f - g\|_{\ell_2(X_N)}$ for $g\in\textrm{span} \{\bar{u}_1, \ldots, \bar{u}_m \}$, $m\leq N$. Moreover, thanks to property \eqref{approx double} of Property \ref{approx discrete properties}, $\bar{R}^m$ can be written in terms of $\ns$-inner products as
\begin{equation}
\bar{R}^{m}(\boldsymbol{x})=\sum_{k=1}^{m} (f, \bar{u}_k)_{\ns} \bar{u}_k(\boldsymbol{x}). 
\label{stab_app}
\end{equation}  
Note that the point \eqref{m=N} of Property \ref{approx discrete properties} proves that $\bar{R}^N = R^N$.

Approximating $R^m$ with its fast computable version $\bar{R}^m$ solves efficiently the problems. We will see in the following sections how to successfully couple this technique with a fast domain decomposition method.
 
\section{Partition of unity method  using stable bases}
\label{PUM_WSVD}
The main idea is to use the stable basis introduced in the  previous section in order to generate local stable approximants  and accumulate them into the global fit.

\subsection{A stable computation of the PU interpolant} \label{PU}

Let $\Omega \subseteq \RR^M$ be an open and bound\-ed domain, and let $\{\Omega_j\}_{j=1}^{d}$ be an open and bounded covering of $\Omega$ satisfying some mild overlap condition among the subdomains (or patches) $\Omega_j$. In other words, the subdomains must form a covering of the domain and, moreover, the overlap must be sufficient so that each interior point $\boldsymbol{x} \in \Omega$ is located in the interior of at least one patch $\Omega_j$.  Furthermore, let us suppose that the set $I(\boldsymbol{x}) = \{j : \boldsymbol{x} \in \Omega_j \}$, for $\boldsymbol{x} \in \Omega$, is uniformly bounded on $\Omega$, with $ \Omega  \subseteq \bigcup_{j=1}^{d} \Omega_j$. Associated with the subdomains we choose partition of unity weight functions $W_j$, i.e. a family of compactly supported, nonnegative and continuous functions subordinate to the subdomain $\Omega_j$, such that $ \sum_{j = 1}^d W_j( \boldsymbol{x})=1$ on $\Omega$ and ${\rm supp}(W_j)  \subseteq \Omega_j$. 

In order to have a better stabilization of the global fit, our goal is to define stable approximants of the form \eqref{stab_app} on each subdomain $\Omega_j$. In other words, for each (local) matrix $A_j \in \mathbb{R}^{N_j \times N_j}$, i.e. the $j$-th interpolation matrix associated with the subdomain $\Omega_j$, a low-rank approximation is computed and thus the so-called WSVD-PU approximant is given by:
\begin{equation} 
{\cal {\bar I}}( \boldsymbol{x})= \sum_{j=1}^{d} {\bar R}_j^{m_j}( \boldsymbol{x} ) W_j ( \boldsymbol{x}), \quad \boldsymbol{x} \in \Omega,
\label{intg}
\end{equation}
where $W_j: \Omega_j \longrightarrow \mathbb{R}$ is a partition of unity weight function and $X_{N_j}= X_N \cap \Omega_j$. As evident, ${\bar R}_j^{m_j}$ defines a local stable RBF approximant on $\Omega_j$ of the form:
\begin{equation}
\bar{R}^{m_j}_j(\boldsymbol{x})=\sum_{k=1}^{m_j} ( \sigma^{(j)}_k)^{-1} (f_{| \Omega_j}, \bar{u}^{(j)}_k)_{\ell_2(X_{N_j})} \bar{u}^{(j)}_k (\boldsymbol{x}), \quad \boldsymbol{x} \in \Omega_j.
\label{loc_fit}
\end{equation}

According to \cite{Wendland02a}, if we assume to have a $ k$-stable partition of unity, then the derivatives of the weight functions satisfy 
\begin{align*}
||D^{ \beta} W_j ||_{L^{ \infty} ( \Omega_j)} \leq \frac{C_{ \beta} }{ \delta_j^{ | \beta|}}, \quad | \beta |	 \leq k, \quad  \forall \beta \in \mathbb{N}^{M},
\end{align*}
where $ \delta_j$ is the diameter of $\Omega_j$ and $C_{ \beta} > 0$ is a constant. As nonnegative functions $W_j \in C^k ( \mathbb{R}^M)$, we may consider a \textsl{Shepard weight}, i.e.,
\begin{align*}
W_j(\boldsymbol{x}) = \frac{\varphi_j(\boldsymbol{x})}{\sum_{k\in I(\boldsymbol{x})} \varphi_k(\boldsymbol{x})}, \quad j=1,\ldots,d,
\end{align*}
$\varphi_j(\boldsymbol{x})$ being compactly supported functions with support on $\Omega_j$, such as the Wendland functions \cite{Wendland05}. 

Before computing the global fit by mean of local stable approximants obtained with the Lanczos procedure, we briefly sketch in the sequel some relevant properties. Since point \eqref{m=N} of Property \ref{approx discrete properties} implies $\bar{R}^{N_j} = R^{N_j}$, we can recover the PU interpolant,  by considering $\bar{R}^{N_j}$, i.e.: 
\begin{equation} 
{\cal I}( \boldsymbol{x})= \sum_{j=1}^{d}  \bar{R}^{N_j}_j( \boldsymbol{x} ) W_j ( \boldsymbol{x}) = \sum_{j=1}^{d} \left( \sum_{k=1}^{N_j} ( \sigma^{(j)}_k)^{-1} (f_{| \Omega_j}, {\bar u}^{(j)}_k)_{\ell_2(X_{N_j})} {\bar u}^{(j)}_k (\boldsymbol{x}) \right) W_j(\boldsymbol{x}) , \quad \boldsymbol{x} \in \Omega.
\label{intg1}
\end{equation}

\begin{remark}
	If the functions $\bar{R}^{N_j}_j$, $j=1,\ldots,d$, satisfy the interpolation conditions $ \bar{R}^{N_j}_j( \boldsymbol{x}_i )= f( \boldsymbol{x}_i)$ for each $\boldsymbol{x}_i \in \Omega_j$, then the global PU approximant inherits the interpolation property of the local interpolants, i.e. 
	\begin{align*}
	{\cal I}( \boldsymbol{x}_i ) = \sum_{j=1}^{d}  \bar{R}^{N_j}_j( \boldsymbol{x}_i ) W_j ( \boldsymbol{x}_i) = \sum_{j \in I(\boldsymbol{x}_i)} f( \boldsymbol{x}_i ) W_j ( \boldsymbol{x}_i) = f( \boldsymbol{x}_i).
	\end{align*}
\end{remark}

{In order to be able to formulate error bounds, we consider the \textsl{fill distance}
	\begin{align*}
	h_{X_N, \Omega} =  \sup_{ \boldsymbol{x} \in \Omega} \min_{ \boldsymbol{x}_i  \in X_N} || \boldsymbol{x} - \boldsymbol{x}_i||_2,
	\end{align*}
and some further assumptions on the regularity of $\Omega_j$.} Specifically, we require that an open and bounded covering $ \{ \Omega_j \}_{j=1}^{d}$ is \textsl{regular} for $( \Omega, X_N)$. This means to fulfill the following properties \cite{Wendland02a}:
\begin{enumerate}
	\item[i.]  for each $ \boldsymbol{ x} \in \Omega$, the number of subdomains $ \Omega_j$ with $ \boldsymbol{x} \in \Omega_j$ is bounded by a global constant $C$;
	\item[ii.]  there exists a constant $C_r > 0$ and an angle $\theta \in (0,\pi/2)$ such that every subdomain $ \Omega_j$ satisfies an interior cone condition with angle $\theta$ and radius $r = C_r h_{ X_N, \Omega}$;
	\item[iii.]  the local fill distances $ h_{ X_{N_j},\Omega_j}$ are uniformly bounded by the global fill distance $h_{X_N, \Omega}$. 
\end{enumerate}

\begin{remark}
	The first property ensures that \eqref{intg1} is actually a sum over at most $C$ summands. Moreover, it is crucial for an efficient evaluation of the global approximant that only a constant number of local interpolants has to be evaluated. It follows that it should be possible to locate those $C$ indices in constant time. The second and third properties are significant for estimating errors of RBF interpolants.
\end{remark}

After defining the space $C_{ \nu}^{k}  ( \mathbb{R}^{M} ) $ of all functions $f \in C^k$ whose derivatives of order $ | \beta |=k $
satisfy $ D^{ \beta} f ( \boldsymbol{x} ) = {\cal O} ( || \boldsymbol{x} ||_2^{ \nu} ) $ for $ || \boldsymbol{x} ||_2 \longrightarrow 0$, we consider the following convergence result {(see Theorem $29.1$ in \cite{Fasshauer07} or also Theorem $15.9$ in \cite{Wendland05})}. It is presented for strictly positive definite functions even if it holds more in general for strictly conditionally positive definite functions.
\begin{theorem}
	Let $ \Omega \subseteq  \mathbb{R}^M$ be open and bounded and suppose that $  X_N= \{ \boldsymbol{x}_i  ,i=1, \ldots, N \} \subseteq \Omega$. Let $ \phi \in C_{ \nu}^{k}  ( \mathbb{R}^{M} ) $ be a
	strictly   positive definite function. Let $ \{ \Omega_j \}_{j=1}^{d}$ be a regular covering for  $( \Omega,  X_N)$  and let $ \{ W_j \}_{j=1}^{d}$ be $k$-stable for $ \{ \Omega_j  \}_{j=1}^{d}$. Then the error between $ f \in {\cal N}_{ \phi} ( \Omega)$, where $ {\cal N}_{ \phi} $ is the native space of  $ \phi $,   and its PU interpolant \eqref{intg1} can be bounded by:
	$$
	| D^{ \beta} f( \boldsymbol{x} ) -  D^{ \beta} {\cal I}( \boldsymbol{x} ) | \leq C h_{  X_N, \Omega}^{\frac{ k+ \nu }{2} - | \beta |} |f|_{{\cal N}_{ \phi} ( \Omega )},
	$$
	for all $ \boldsymbol{x} \in \Omega $ and all $ | \beta | \leq k/2$.
\end{theorem}

\subsection{The PU algorithm: the Lanczos procedure} \label{lancz}
For each subdomain, in order to generate the local stable approximation matrix, the Lanczos method is applied to the matrix $A_j \in \mathbb{R}^{N_j \times N_j}$ and to the function values $\boldsymbol{f}_j \in \mathbb{R}^{N_j}$ associated with the subdomain $\Omega_j$.  In this way the matrix $H_{m_j}$ and the Lanczos basis $\{ \boldsymbol{p}^{(j)}_i \}^{m_j}_{i=1}$ are computed for each subdomain. Then, for each  interpolation problem a  local stable basis is formed.

By using a different stopping criterion in the Lanczos algorithm, with respect to the one employed in \cite{Demarchi15}, we can compute stable bases for a wider family of RBFs, both globally defined and compactly supported. The main problem in the Lanczos procedure concerns the stopping criterion, (see  {\fontfamily{pcr} \selectfont Step 3} of the {\tt Lanczos Algorithm}). From Property \ref{discrete properties} (point \eqref{discrete sum condition}) and Property \ref{approx discrete properties} (point \eqref{m=N}) a reliable  one is:
\begin{equation}
\bigg|\phi(0)-\dfrac{1}{N_j} \sum_{k=1}^{m_j} \alpha^{(j)}_k \bigg|< \tau,
\label{toll}
\end{equation}
for a certain fixed tolerance $\tau$, which is supposed to be equal  for all the subdomains. We point out that, even if $\tau$ is fixed among the subdomains, the left-hand side of \eqref{toll} depends on the specific patch. Moreover, supposing to have a quasi-uniform node distribution, there are no restrictions in selecting the same tolerance for all the subdomains. On the opposite, if the points are more clustered in several subdomains, one can always keep a fixed tolerance, but techniques enabling us to select suitable  sizes of the different subdomains are recommended (see e.g. \cite{Safdari}).

From Property \ref{approx discrete properties} (point \eqref{m=N}), the fact that we impose as maximum number of iterations, in the {\tt Lanczos} {\tt Algorit-}\\{\tt hm} at {\fontfamily{pcr} \selectfont Step 2}, exactly the number of nodes in $\Omega_j$, i.e. $N_j$, naturally follows. 

\begin{table}[!ht]
	{\normalsize
		\begin{center}
			\begin{tabular}{p{12cm}*{1}{c}}
				\hline
				\vskip 0.01 cm 
				{\fontfamily{pcr} \selectfont INPUTS:} $N_j$, number of data in $\Omega_j$;  $A_j$, the local interpolation matrix;
				\vskip 0.08cm
				\hskip 2.2 cm $\boldsymbol{f}_j$, the function values associated to  $\Omega_j$;
				\vskip 0.08cm
				\hskip 2.2 cm $\tau$, the tolerance used as stopping criterion; $\phi$, the radial basis
				\vskip 0.08cm
				\hskip 2.2 cm  function. 
				\vskip 0.08cm
				{\fontfamily{pcr} \selectfont OUTPUTS:} $\boldsymbol{p}^{(j)}_1, \ldots,\boldsymbol{p}^{(j)}_m$, the new basis in $\Omega_j$; $H_{m_j}$, the tridiagonal matrix.
				\vskip 0.12 cm
				{\fontfamily{pcr} \selectfont Step 1:}
				Set $\beta^{(j)}_1=0$; $\boldsymbol{p}^{(j)}_0=0$; $\boldsymbol{p}^{(j)}_1=\dfrac{\boldsymbol{f}_j}{|| \boldsymbol{f}_j||_2}$.
				\vskip 0.12 cm
				{\fontfamily{pcr} \selectfont Step 2:}
				\textbf{For} $i=1:N_j$
				\vskip 0.08cm
				\hskip 2.2 cm $\boldsymbol{\tilde{p}}^{(j)}_i=A_j\boldsymbol{p}^{(j)}_i-\beta_i \boldsymbol{p}^{(j)}_{i-1}$
				\vskip 0.08cm
				\hskip 2.2 cm $\alpha^{(j)}_i= ( \boldsymbol{\tilde{p}}^{(j)}_i,\boldsymbol{p}^{(j)}_i)$
				\vskip 0.08cm
				\hskip 2.2 cm $\boldsymbol{\tilde{p}}^{(j)}_i=\boldsymbol{\tilde{p}}^{(j)}_i-\alpha^{(j)}_i \boldsymbol{p}^{(j)}_{i}$
				\vskip 0.08cm
				\hskip 2.2 cm $\beta^{(j)}_{i+1}=||\boldsymbol{\tilde{p}}^{(j)}_i||_2$
				\vskip 0.12 cm
				\hskip 2.2 cm {\fontfamily{pcr} \selectfont Step 3:}
				\textbf{If} $\beta^{(j)}_{i+1}=0$ or $|\phi(0)-\frac{1}{N_j} \sum_{k=1}^{i} \alpha^{(j)}_k|< \tau$
				\vskip 0.08cm
				\hskip 4.2 cm \textbf{break}
				\vskip 0.12 cm
				\hskip 2.2 cm $\boldsymbol{p}^{(j)}_{i+1}=\boldsymbol{\tilde{p}}^{(j)}_i/ \beta^{(j)}_{i+1}$
				\\[\smallskipamount]
				\hline
			\end{tabular}
		\end{center}
		\caption{The {\tt Lanczos Algorithm}. Routine performing the Lanczos 
			procedure.}
		\label{Lan_code}
	}
\end{table}

Then, once the matrix $H_{m_j}$ is found for $\Omega_j$ the stable basis is computed by calculating the singular value decomposition of $H_{m_j}$ and a local approximant on each subdomain in the form \eqref{loc_fit} is computed. Then the local fits are accumulated into the global one.

The use of stable bases by decomposing the initial problem into many small ones leads to a larger benefit in terms of accuracy than employing a global approach. In fact if one uses a global method  the approximant results stable, but a large number of terms in the Lanczos procedure are neglected. This surely leads to a decrease of the fit accuracy. Whereas the local method turns out to be really accurate since, dealing with small problems, less terms in the computation of the basis are eliminated to preserve stability. 

Extensive numerical experiments support our findings.

\subsection{The PU algorithm: the block-based algorithm structure} \label{PU_ALG}

The key step of the PU method consists in organizing the scattered data among the subdomains. To this aim the kd-tree partitioning structures are widely used, \cite{Fasshauer07}. However, they are not specifically implemented for the PU method.

Here a novel partitioning procedure, specifically the so-called \emph{block-based partitioning structure},  built for bivariate and trivariate interpolation in order to determine the points belonging to the different PU subdomains, is considered. Even if such partitioning structure is robust enough to work on 2D or 3D irregular domains, \cite{Cavoretto15c},  we present such efficient technique for a scattered data set lying in the unit square, i.e. $\Omega=[0,1]^2$. 

At first, a partition of unity structure, composed by $d$ circular patches $\Omega_j$ of radius:
\begin{equation} 
\delta = \sqrt{ \dfrac{2}{d}},
\label{PU_radius}
\end{equation}
and whose centres $\bar{\boldsymbol{x}}_j  ,j=1, \ldots, d$, are a grid of points on $\Omega$, is generated.
As in \cite{Fasshauer07}, the number of PU subdomains is chosen so that $N/d \approx 4$. This choice and  \eqref{PU_radius} lead to a reliable partition of unity structure since, in this way,  patches form a covering  of the domain $\Omega$.  

In order to find the points belonging to the different subdomains and consequently solve, with the use of stable bases, $d$ small interpolation problems,  we propose a new partitioning structure. It leads to a natural searching procedure that turns out to be really cheap in terms of computational complexity. To this aim we first cover $\Omega$ with $q^2$ square blocks, where the number $q$ of blocks along one side of the unit square is:
\begin{equation} 
q= \bigg \lceil \frac{\displaystyle 1}{\displaystyle \delta} \bigg \rceil.
\label{q}
\end{equation}
In this way the width of blocks is equal to the subdomain radius. This choice can appear trivial, but on the contrary it enables us to consider in the searching process an optimized number of blocks. 

Blocks are numbered from $1$ to $q^2$ (bottom to top, left to right). Thus, with a repeated use of a quicksort routine  the set $ X_N$ is partitioned by the block-based partitioning structure into $q^2$ subsets $X_{N_k}$, $k = 1, \ldots , q^2$, where $X_{N_k}$  are the points stored  in the $k$-th  \emph{neighbourhood}, i.e. in the $k$-th block and in its eight neighbouring blocks. In such framework, we are able to get an optimal procedure to find the nearest points. In fact, given a subdomain $\Omega_j$, whose centre belongs to the $k$-th block, we search for all data lying in the  $j$-th subdomain only among those lying in  the $k$-th neighbourhood.

\begin{remark}
	The same partitioning structure, in case of Compactly Supported RBFs (CSRBFs), must be considered locally for each subdomain. In fact, in order to build the $j$-th stable approximation matrix, among all points lying in the $j$-th subdomain, only those belonging to the support of the CSRBF must be considered.
\end{remark}

\begin{remark}
	Among several routines which can be employed to determine the neighboring points, we choose the block-based data structure. Anyway, we stress that the algorithm, here proposed, works in any dimension $M$, while the block-based data structure is only implemented for $M=2,3$, \cite{Cavoretto15c}. Thus in higher dimensions such structure must be replaced by standard routines, such as kd-trees, \cite{Arya98,deBerg97,Fasshauer07}.
\end{remark}

\begin{remark}
		For easiness of the reader, the procedure is here described in the unit square, however, following \cite{Cavoretto15c}, any extension to irregular domains is possible. 
\end{remark}

\section{Complexity analysis}\label{compl_anal_WSVD}



Since the stable WSVD-PU algorithm is characterized by the construction of local RBF stable approximants, we consider the local data sets, composed by $N_j$ points, $j = 1, 2, \ldots , d$. Thus, the complexity of this algorithm is influenced by the following computational issues:
\begin{enumerate}
	\item[i.] organize by means of a partitioning structure the nodes among the subdomains,
	\item[ii.] compute the stable basis on each subdomain
\end{enumerate}
Concerning the efficient organization of points, an extensive complexity analysis,  briefly shacked in  Subsection \ref{comp_PS}, can be found in \cite{Cavoretto15c}. The cost associated to the computation of a local stable basis is investigated in Subsection \ref{compl_lan}.

\subsection{Computation of a stable basis}
\label{compl_lan}
Performing the Lanczos procedure on a matrix $B \in \mathbb{R}^{n \times n}$ requires ${\cal O} (kn^2)$, where $k$ is the number of vectors computed by the algorithm, i.e. $k$ is the \emph{good} low rank approximation, (a priori unknown in our case), \cite{Chen09}.

Given $A_j \in \mathbb{R}^{N_j \times N_j}$ the interpolation matrix defined on the $\Omega_j$, the Lanczos method forms the matrix $H_{m_j}$ for  $\Omega_j$ after $m_j$ iterations.
Usually we have $m_{j} \ll N_{j}$, but in some cases the maximum number of iterations $N_j$ can be reached and so, in a more general setting, $m_{j}\leq N_{j}$. This routine requires:
\begin{equation}
{\cal O} (m_j N^2_j) \leq {\cal O} (N^3_j),
\label{comp_lan}
\end{equation}
time complexity.
Thus for each subdomain  the upper bound for the computational time of the Lanczos procedure is given by the right-hand side of \eqref{comp_lan}. 

In case of sparse matrices, such as the  ones arising from the use of  CSRBFs, the  Lanczos procedure can be performed in:
${\cal O} (m_j(N_j+ \tilde{n}))$ time complexity, where $\tilde{n}$  is the number  of non-zero entries. 

Then a singular value decomposition is applied to the matrix $H_{m_j}$. We remark that performing a singular value decomposition on a matrix  $B \in \mathbb{R}^{n \times k}$ requires ${\cal O} (4n^2k+8nk^2+9k^3)$ time complexity.

The singular value decomposition  for each subdomain is applied to the matrix $H_{m_j}$; once more we stress that $m_j \ll N_j$. Thus for each subdomain the singular value decomposition can be performed in:  
\begin{equation}
{\cal O} (4m^2_jm_j+8m_jm^2_j+9m^3_j) \approx {\cal O} (m^3_j) 
\label{comp_svd}
\end{equation}
time complexity.

\subsection{The partitioning structure}
\label{comp_PS}
Let us now focus on the block-based partitioning structure used to organize the $N$  data sites in blocks. We remark that such efficient organization of points is specifically implemented for 2D data sets. Anyway, the proposed WSVD-PU algorithm is robust enough to work in any dimension $M$, provided that a different partitioning structure is performed.

Let ${\bar n}_s$  be  the number of data sites belonging to a strip. The procedure used to store the points among the different subdomains is based on recursive calls to a \emph{quicksort} routine which requires $ {\cal O} ({\bar n} \log {\bar n})$, where ${\bar n}$ is the number of elements to be sorted. Thus, letting $N/q$ the average number of points lying in a strip,  the computational cost needed to organize the $N$ points among the different subdomains is:
\begin{equation}
{\cal O} \bigg( N \log{N} + \sum_{s=1}^{q} {\bar n}_s \log{{\bar n}_s} \bigg) \approx {\cal O} \bigg( \frac{\displaystyle  3}{\displaystyle  2} N \log{ N} \bigg).
\label{CA1}
\end{equation}


Concerning the searching procedure, for each subdomain a quicksort procedure is used to order distances. Thus observing that the data sites  in a neighbourhood are about $N/(3q)^2$ and taking into account the definitions of $q$ and $\delta$, the complexity can be estimated by:
\begin{equation}
{\cal O} \bigg( \frac{ \displaystyle N}{\displaystyle (3q)^2} \log{ \frac{ \displaystyle N}{\displaystyle (3q)^2}}\bigg) \approx {\cal O} \bigg( \frac{ \displaystyle 2 N}{\displaystyle 9d} \log{ \frac{ \displaystyle 2N}{\displaystyle 9d}}\bigg) \approx {\cal O} (1).
\label{CA7}
\end{equation}
The estimate \eqref{CA7} follows from the fact that we built a partitioning structure strictly related to the size of the subdomains and ad hoc for the PU method. 
\begin{remark}
	The same computational cost \eqref{CA7}, in case of CSRBFs, must be considered locally for each subdomain,  to build the sparse interpolation and evaluation matrices. In such steps we usually have a relatively small number of nodes $N_j$, with $N_j \ll N$,   where the index $j$ identifies the $j$-th subdomain.
\end{remark}
\section{Numerical experiments} \label{NE_LAN}
This section is devoted to point out, by means of extensive numerical simulations,  stability and accuracy of the WSVD-PU interpolant. To this aim comparisons with the standard PU interpolant will be carried out.

Experiments are performed considering $N=(2^k + 1)^2,$ $k = 6, 7, 8$, uniformly random Halton nodes, a grid of $d = \lfloor \sqrt{N} /2 \rfloor^2$ subdomain centres and a grid of $s = 40 \times 40$ evaluation  points, which are contained in the unit square $\Omega = [0, 1] \times [0, 1]$.

In order to show the high stability of the proposed method, we compute the Root Mean Square Error (RMSE), i.e.
\begin{eqnarray} \label{RMSE}
	\textrm{RMSE} = \sqrt{\frac{1}{s}\sum_{i=1}^{s} |f(\tilde{\boldsymbol{x}}_i) - {\bar {\cal I}}(\tilde{\boldsymbol{x}}_i)|^2},
\end{eqnarray}
for different values of the shape parameter in the range $\varepsilon = [10^{-3}, 10^{2}]$. Moreover, in order to point out the versatility of the proposed method,  different kernels with different order of smoothness are considered, see Table \ref{tab_1}. The error \eqref{RMSE} is computed using as test function the well-known Franke's function:
\begin{align*}
f(x_1,x_2)&=\frac{3}{4}\exp\left[-\frac{(9x_1-2)^2+(9x_2-2)^2}{4}\right]+\frac{3}{4}\exp\left[-\frac{(9x_1+1)^2}{49}-\frac{9x_2+1}{10}\right]\\
&+\frac{1}{2} \exp\left[-\frac{(9x_1-7)^2+(9x_2-3)^2}{4}\right]-\frac{1}{5} \exp\left[-(9x_1-4)^2-(9x_2-7)^2\right].
\end{align*}

\begin{figure}[ht!]
\begin{center}
\includegraphics[width=8.8cm]{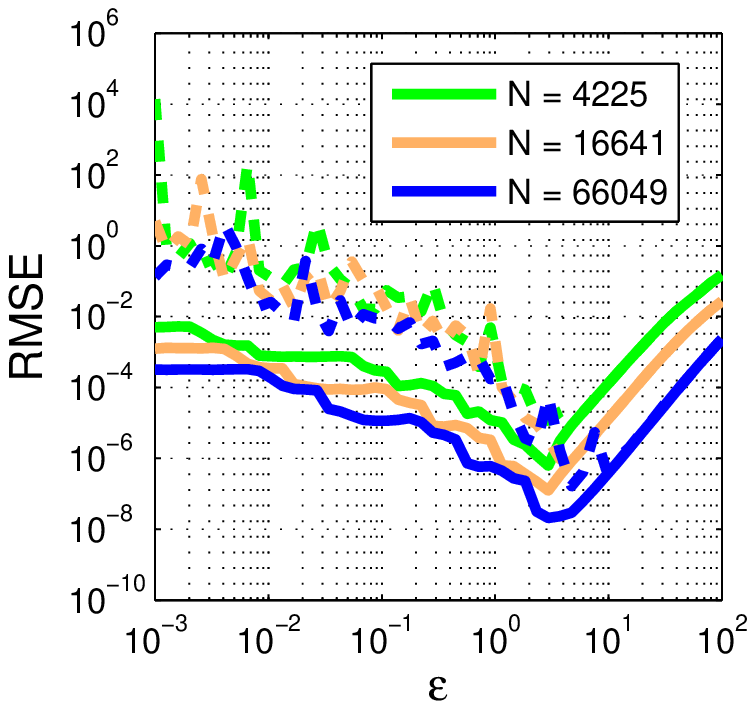} \hskip -1.6cm
\includegraphics[width=8.8cm]{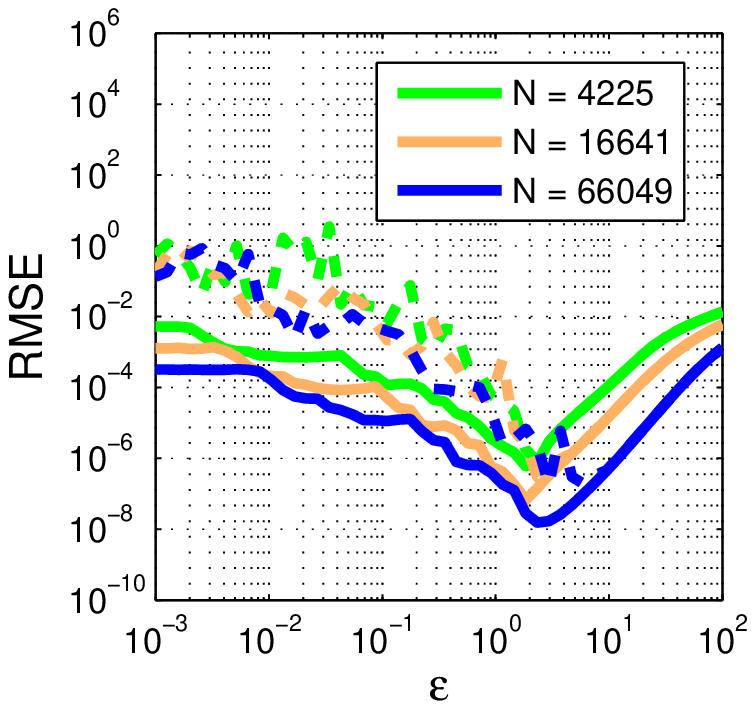} \\
\includegraphics[width=8.8cm]{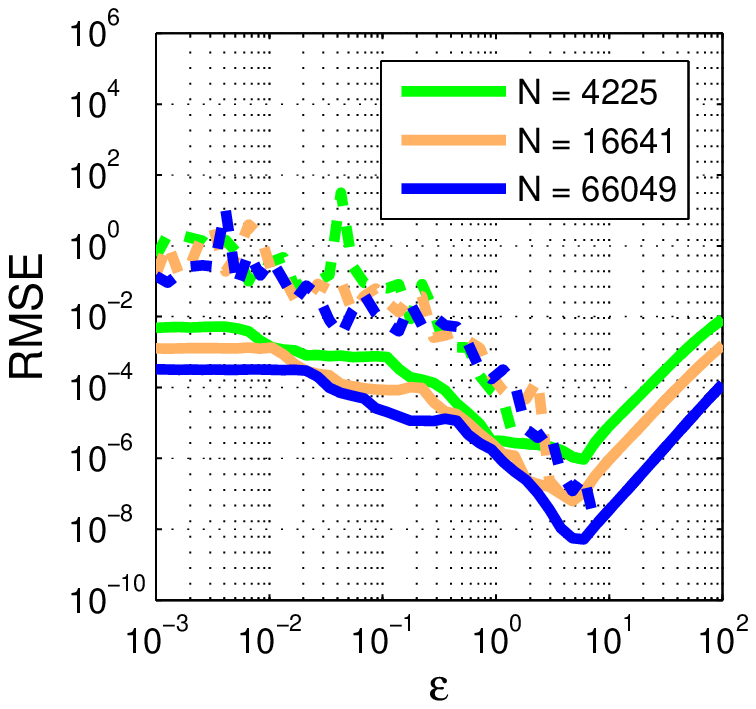} \hskip -1.6cm
\includegraphics[width=8.8cm]{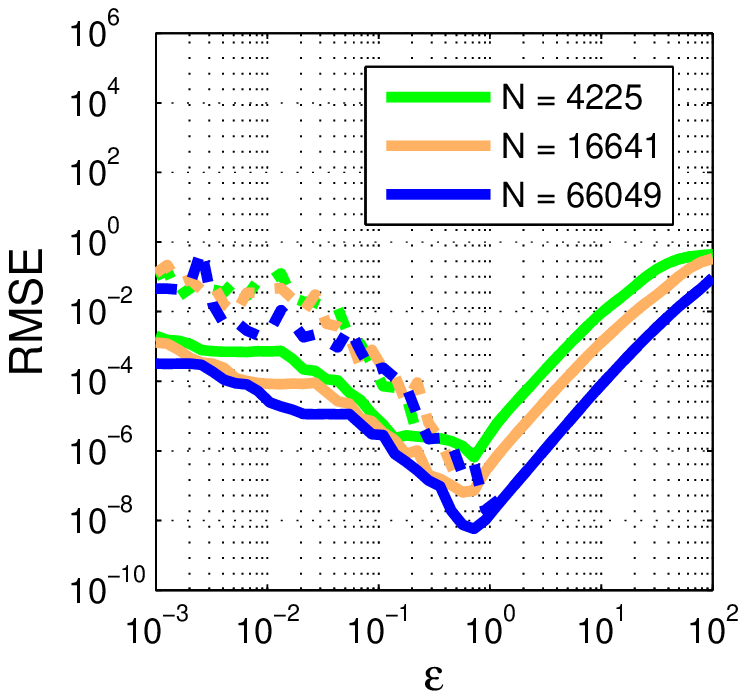} \\
\includegraphics[width=8.8cm]{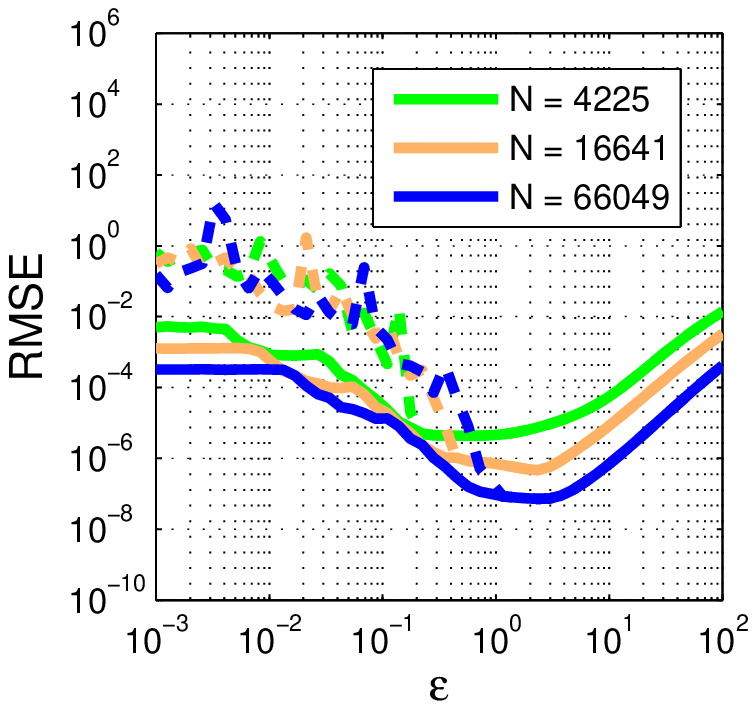} \hskip -1.6cm
	\includegraphics[width=8.8cm]{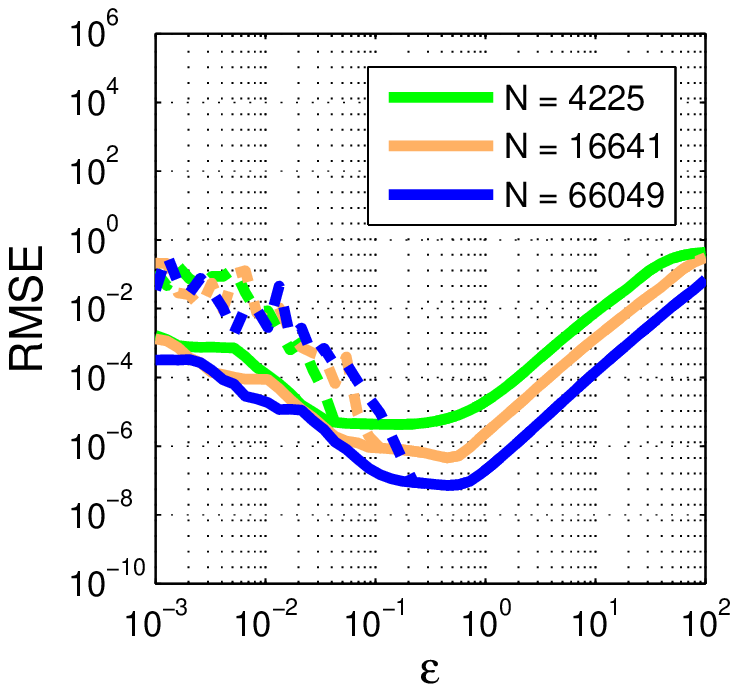}
\end{center}
\caption{RMSEs obtained by varying $\varepsilon$ for $C^\infty$, $C^6$ and $C^4$ kernels. The classical PU interpolant is plotted with dashed line and the  WSVD-PU approximant with solid line. {From left to right, top to bottom we consider the GA, IMQ, M6, W6, M4 and W4 kernels, respectively.}}
\label{fig_1}
\end{figure}

In Figure \ref{fig_1}  we compare the RMSEs obtained by means of the WSVD-PU interpolant (solid line) with the ones obtained performing the classical PU method (dashed line).  As tolerance value in \eqref{toll} we set $10^{-14}$.
These graphs point out that the use of the WSVD-PU local approach reveals a larger stability than the standard PU interpolant.
Moreover, the use of a local method enables us to improve the RMSE for the optimal shape parameter in case of flat kernels, see Figure \ref{fig_1} and Table \ref{tab_1}. This is consistent with the fact that in a local stable method, differently from \cite{Demarchi15}, we have to solve small linear systems and therefore few terms are neglected in \eqref{loc_fit}.
Furthermore, from Figure \ref{fig_1}   we can note that the WSVD-PU method turns out to be more effective with flat kernels, while for more picked bases the improvement of using stable bases becomes negligible as the order of bases function decreases.
Thus, from our numerical experiments, we can observe three kinds of behavior depending on different RBF regularity classes. Specifically, the features of such classes, which differ both in terms of stability and accuracy from the standard basis, can be summarized as:
\begin{itemize}
	\item[i.] for $C^{\infty}$ kernels: improvement of stability and of the optimal accuracy;
	\item[ii.] for $C^k$ kernels, with $k\geq 1$: improvement of stability and same optimal accuracy;
	\item[iii.] for $C^0$ kernels: same stability and same optimal accuracy.
\end{itemize}

\begin{table}[ht!]
		\begin{center}
			\begin{tabular}{|c|c|c|c|c|c|c|c|c|c|} \hline
			  & & \multicolumn{2}{c|}{  \rule[-2mm]{0mm}{7mm}  GA} & \multicolumn{2}{c|}{  \rule[-2mm]{0mm}{7mm}  IMQ} & \multicolumn{2}{c|}{  \rule[-2mm]{0mm}{7mm} M6} & \multicolumn{2}{c|}{  \rule[-2mm]{0mm}{7mm} W6} \\
			  \cline{3-10} \rule[-2mm]{0mm}{7mm}
			  $N$ & method	& RMSE & $\varepsilon_{opt}$ & RMSE & $\varepsilon_{opt}$ & RMSE & $\varepsilon_{opt}$ & RMSE & $\varepsilon_{opt}$ \\
				\hline 
				\rule[0mm]{0mm}{3ex}
$\hskip-2pt 4225$ & PU & $1.16{\rm E}-5$  &   $2.95$    & $8.20{\rm E}-7$  &   $2.33$	   & $9.34{\rm E}-7$  &   $5.96$     & $6.64{\rm E}-7$  &   $0.72$   \\
                  &  WSVD-PU            & $6.20{\rm E}-7$  &   $2.95$    & $5.98{\rm E}-7$  &   $1.84$	   & $9.34{\rm E}-7$  &   $5.96$     & $6.64{\rm E}-7$  &   $0.72$   \\
        \hline
				\rule[0mm]{0mm}{3ex}
$16641$     &  PU  & $9.70{\rm E}-7$  &   $3.73$    & $2.94{\rm E}-7$  &   $2.33$	   & $6.18{\rm E}-8$  &   $4.71$     & $6.44{\rm E}-8$  &   $0.57$   \\
            &  WSVD-PU  & $1.25{\rm E}-7$  &   $2.95$    & $6.78{\rm E}-8$  &   $1.84$	   & $6.20{\rm E}-8$  &   $4.71$     & $6.49{\rm E}-8$  &   $0.57$   \\
        \hline
			  \rule[0mm]{0mm}{3ex}
$66049$     &  PU   & $1.64{\rm E}-7$  &   $4.71$    & $1.78{\rm E}-7$  &   $2.94$	   & $1.28{\rm E}-8$  &   $7.54$     & $2.03{\rm E}-8$  &   $0.91$   \\
            &  WSVD-PU   & $2.09{\rm E}-8$  &   $2.95$    & $1.54{\rm E}-8$  &   $2.33$	   & $5.10{\rm E}-9$  &   $5.96$     & $5.70{\rm E}-9$  &   $0.72$   \\
			  \hline 
			\end{tabular}
		\end{center}
			\caption{RMSEs obtained by using optimal values of $\varepsilon$ for {GA, IMQ, M6 and W6} kernels.}
			\label{tab_1}
	\end{table}

Moreover, since we  are interested in pointing out the efficiency of the proposed WSVD-PU algorithm, in Table \ref{tab_3} we also report the CPU times obtained by using our stable interpolation method with the Gaussian RBF as local approximant, for each of the three different data sets. Tests have been carried out on a Intel(R) Core(TM) i3 CPU M330 2.13 GHz processor.

\begin{table}[ht!]
		\begin{center}
			\begin{tabular}{|c|c|c|c|} \hline
                \rule[0mm]{0mm}{3ex}
			  $N$ & $4225$	& $16642$ & $66049$  \\
				\hline 
				\rule[0mm]{0mm}{3ex}
$\hskip-2pt $ CPU [s] & $3.94$ & $14.69$  &   $57.73$   \\
			  \hline 
			\end{tabular}
		\end{center}
			\caption{CPU times (in seconds) for the PU-WSVD method.}
			\label{tab_3}
	\end{table}	
	
\section{Application to real world data sets}
\label{applicazione}

\begin{figure}[ht!]
	\begin{center}
			\includegraphics[width=8.8cm]{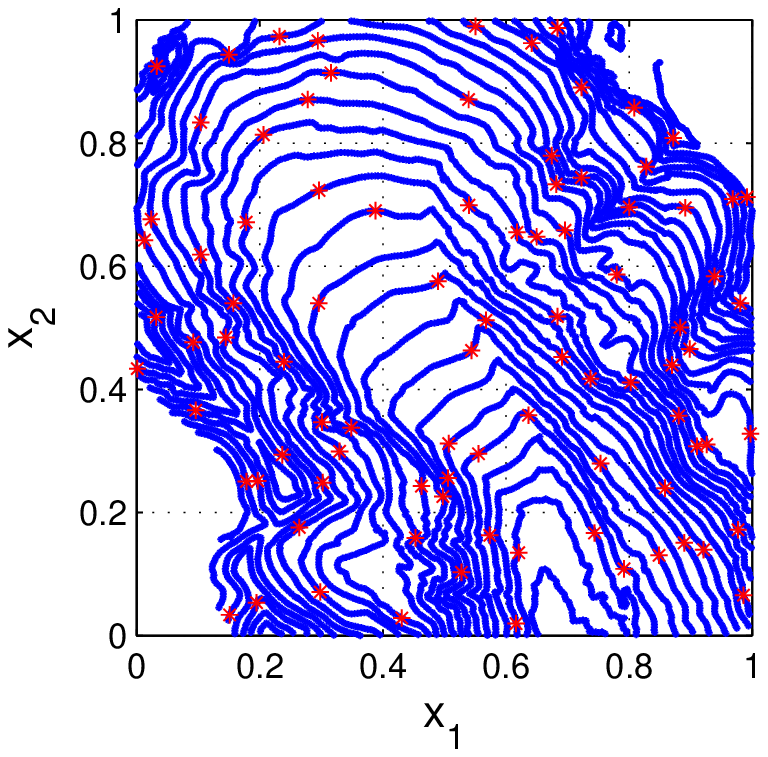} \hskip -1.6cm
			\includegraphics[width=8.8cm]{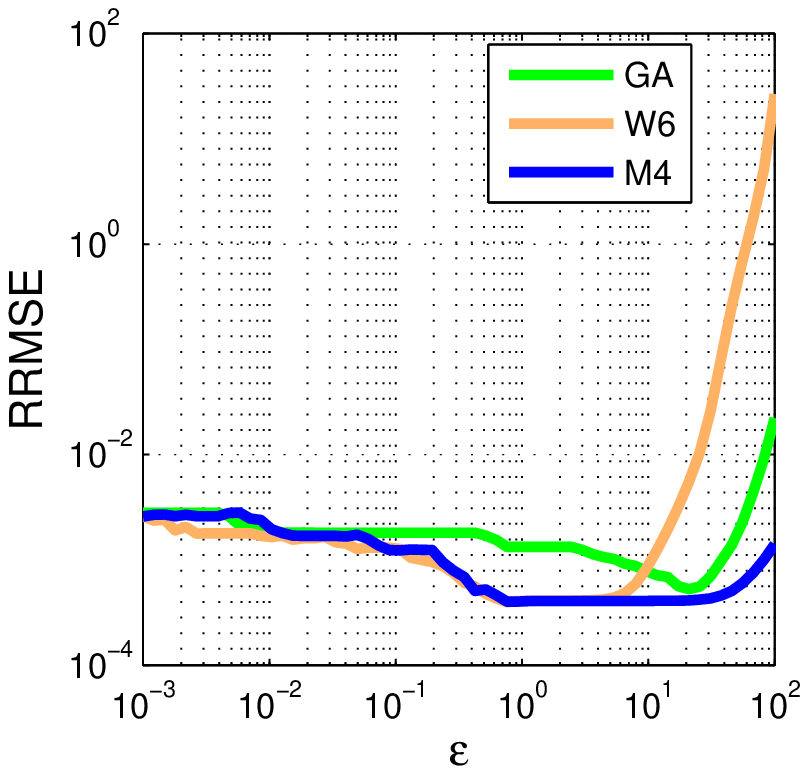}
	\end{center}
	\caption{Left: a 2D view of the glacier data set. Right: RRMSEs obtained by varying $\varepsilon$ with {GA, W6 and M4} kernels for the  WSVD-PU approximant.}
	\label{fig_2}
\end{figure}

In this section we focus on an application to Earth's topography, which consists in approximating with our algorithm a set of real scattered data. In particular, we consider  the so-called \emph{glacier}  data set. It is composed by $8345$ points  representing  digitized height contours of a glacier, \cite{davydov,Wendland01}. The difference between the highest and the lowest point  is $800$ m. Such points, differently from the Halton data, are not quasi-uniform. Furthermore, they are distributed on an irregular domain $\Omega \subseteq [0,1]^2$. A 2D view of such points is plotted in Figure \ref{fig_2} (left).  Among them we select $s=90$ points for the cross-validation (plotted in red in the left frame of Figure \ref{fig_2}).

Since the nodes of the glacier data set do not cover the whole unit square, after generating an initial grid of subdomain centres in $[0,1]^2$, we reduce such points taking only those lying in $\Omega$ by means of the technique described in \cite{Cavoretto15c}.

To obtain reliable and numerically significant results of the error, in this application it is more appropriate to use  the Relative RMSE (RRMSE):
\begin{align*} 
\textrm{RRMSE} = \sqrt{\dfrac{1}{s}\sum_{i=1}^{s} \dfrac{\left|f(\tilde{\boldsymbol{x}}_i) - {\bar {\cal I}}(\tilde{\boldsymbol{x}}_i)\right|^2}{\left|f(\tilde{\boldsymbol{x}}_i)\right|^2}}.
\end{align*}

In Figure \ref{fig_2}, we show how the RRMSEs vary with respect to the shape parameter $\varepsilon \in [10^{-3}, 10^{2}]$. In doing so, we consider {the following kernels: GA, W6 and M4}.

As already shown, the results point out once more that the proposed approach is stable  and moreover turns to be effective also in  applications. The errors for the optimal shape parameter $\varepsilon$ are shown in Table \ref{tab_2}. Since we refer to points with highly varying densities and thus truly ill-conditioned matrices, the classical PU method does not give acceptable approximations. Consequently, we do not report the errors  obtained with this standard algorithm.

\begin{table}[ht!]
	\begin{center}
		\begin{tabular}{|c|c|c|c|c|c|} \hline
			\multicolumn{2}{|c|}{  \rule[-2mm]{0mm}{7mm}  GA} & \multicolumn{2}{c|}{  \rule[-2mm]{0mm}{7mm}  W6} & \multicolumn{2}{c|}{  \rule[-2mm]{0mm}{7mm} M4} \\
			\cline{1-6} \rule[-2mm]{0mm}{7mm}
			RRMSE & $\varepsilon_{opt}$ & RRMSE & $\varepsilon_{opt}$ & RRMSE & $\varepsilon_{opt}$  \\
			\hline 
			\rule[0mm]{0mm}{3ex}
			 $5.26{\rm E}-4$  & $20.9$  &   $3.96{\rm E}-4$  & $0.76$   &   $4.02{\rm E}-4$  & $0.76$   \\
			\hline 
		\end{tabular}
	\end{center}
	\caption{RRMSEs obtained by using optimal values of $\varepsilon$ for {GA, W6 and M4} kernels.}
	\label{tab_2}
\end{table}




\section*{Acknowledgements}  
{The first,} third and fourth authors are partially supported by the University of Torino through research project \lq\lq Metodi numerici nelle scienze applicate\rq\rq. The second and fifth authors are partially supported by the funds of the University of Padova, project CPDA124755 \lq\lq Multivariate approximation with application to image reconstruction\rq\rq.



\end{document}